\documentclass[10pt,a4paper, leqno]{article}
\usepackage[utf8]{inputenc}
\usepackage{amsmath}
\usepackage{amsfonts}
\usepackage{amssymb}
\usepackage{comment}
\usepackage{hyperref}
\usepackage{enumitem,pstricks,xy,pst-node}
\xyoption{all}
\usepackage{amsthm}
\usepackage{booktabs} %to produce nicer tables
\usepackage[color]{changebar} %to introduce changebars
\cbcolor{red} %the color of the changebars
\usepackage{mathrsfs} %to use \mathscr
\usepackage{verbatim} %to use \verb
\usepackage[disable] %to do notes in the document
{todonotes}
\usepackage[normalem]{ulem} %to strike out text
\usepackage{caption} %to put some space between tables and caption
\usepackage{nicefrac} %to use nice inline fractions, e.g. 1/2
\usepackage{enumitem} % to change indent of itemize

\newtheorem{theorem}{Theorem}[section]

\newtheorem{definition}[theorem]{Definition}
\theoremstyle{definition}

\newtheorem{remark}[theorem]{Remark}
\newtheorem*{ack}{Acknowledgements}

\theoremstyle{remark}

\newcommand{\cH}{\mathcal{H}}
\newcommand{\cW}{\mathcal{W}}
\newcommand{\Gr}{\mathbf{Gr}}
\newcommand{\PP}{\mathbb{P}}

\DeclareMathOperator{\Spec}{Spec}
\DeclareMathOperator{\Hb}{H}
\DeclareMathOperator{\hsm}{h}
\DeclareMathOperator{\pgl}{PGL}
\DeclareMathOperator{\Pic}{Pic}

\usepackage[textwidth=125mm,textheight=185mm]{geometry}

\newcommand{\ladi}{\begin{lastadd}}
\newcommand{\ladf}{\end{lastadd}}
\newcommand{\lrei}{\begin{lastrem}}
\newcommand{\lref}{\end{lastrem}}
\newenvironment{lastadd}
{\cbstart\color{red}}
{\todo{red to remove}\cbend}
\newenvironment{lastrem}
{\cbstart\color{yellow}}
{\cbend}

\author{Hanieh Keneshlou\thanks{Mathematik und Informatik, Universit\"at des Saarlandes, Campus E2 4, D-66123 Saarbr\"ucken, Germany. keneshlou@math.uni-sb.de} \and Fabio Tanturri\thanks{Corresponding author. Laboratoire Paul Painlev\'e, Universit\'e de Lille, 59655 Villeneuve d'Ascq CEDEX, France. Fabio.Tanturri@math.univ-lille1.fr}}
\title{The unirationality of the Hurwitz schemes $\mathcal{H}_{10,8}$ and $\mathcal{H}_{13,7}$}
\begin{document}
\maketitle

\begin{abstract}
\noindent We show that the Hurwitz scheme $\mathcal{H}_{g,d}$ parametrizing $d$-sheeted simply branched covers of the projective line by smooth curves of genus $g$, up to isomorphism, is unirational for $(g,d)=(10,8)$ and $(13,7)$. The unirationality is settled by using liaison constructions in $\mathbb{P}^1 \times \mathbb{P}^2$ and $\mathbb{P}^6$ respectively, and through the explicit computation of single examples over a finite field.
\end{abstract}

\section*{Introduction}
The study of the birational geometry of the moduli spaces of curves together with additional data such as marked points or line bundles is a central subject in modern algebraic geometry. For instance, understanding the geometry of the Hurwitz schemes
\[
\mathcal{H}_{g,d}:=\{\xymatrix{C \ar[r]^-{d:1} & \mathbb{P}^1} \mbox{ simply branched cover } |\; C \mbox{ smooth of genus }g\}/\sim
\]
parametrizing $d$-sheeted simply branched covers of the projective line by smooth curves of genus $g$, up to isomorphism, has an important role in shedding light on the geometry of the moduli spaces of curves $\mathcal{M}_g$. It was through Hurwitz spaces that Riemann \cite{Riemann57} computed the dimension of $\mathcal{M}_g$, and Severi \cite{Severi68}, building on works of Clebsch and L\"uroth \cite{Clebsch73}, and Hurwitz \cite{Hurwitz91}, gave a first proof that $\mathcal{M}_g$ is irreducible.

Recently, the birational geometry of Hurwitz schemes has gained increasing interest, especially concerning their unirationality. By classical results of Petri \cite{Petri}, Segre \cite{Segre}, and Arbarello and Cornalba \cite{ArbarelloCornalba}, it has been known for a long time that $\mathcal{H}_{g,d}$ is unirational in the range $2 \leq d \leq 5$ and $g \geq 2$. For $g \leq 9$ and $d\geq g$, the unirationality has been proved by Mukai \cite{Mukai1995}. The most recent contributions have been given by \cite{VerraUnirationality, GeissUnirationality, GeissThesis, SchreyerComputer, SchreyerTanturriMatrix, DamadiSchreyer} and show how active this research area is. For a more complete picture on the unirationality of Hurwitz spaces, the related speculations and open questions, we refer to \cite{SchreyerTanturriMatrix}. Hurwitz spaces and their Kodaira dimension are also considered in the very recent \cite{FarkasEffectiveDiv}.

The main contribution of this paper is the proof of the unirationality of the Hurwitz schemes $\mathcal{H}_{10,8}$ and $\mathcal{H}_{13,7}$ (Theorem \ref{H108unirational} and Theorem \ref{H137unirational}). In \cite[\textsection 1]{SchreyerTanturriMatrix} it is speculated that $\mathcal{H}_{g,d}$ is unirational for pairs $(g,d)$ lying in a certain range: we remark that our two cases lie in that range, and respect perfectly this speculation.

The key ingredient for both results is the construction of dominant rational families of curves constructed via liaison in $\mathbb{P}^1 \times \mathbb{P}^2$ and $\mathbb{P}^6$ respectively. The proof of the unirationality of $\mathcal{H}_{10,8}$ is based on the observation that a general $8$-gonal curve of genus 10 admits a model in $\mathbb{P}^1 \times \mathbb{P}^2$ of bidegree $(6,10)$, which can be linked in two steps to the union of a rational curve and five lines. We show that this process can be reversed and yields a unirational parametrization of $\mathcal{H}_{10,8}$.

For $\mathcal{H}_{13,7}$, we use the fact that a general $7$-gonal curve of genus $13$ can be embedded in $\mathbb{P}^6$ as a curve of degree $17$, which is linked to a curve $D$ of genus $10$ and degree $15$. We show that also this process can be reversed; to exhibit a unirational parametrization of such $D$'s, we prove the unirationality of $\mathcal{M}_{10,n}$ for $n \leq 5$ (Theorem \ref{M10n}), a result of independent interest, and we use a general curve together with $3$ marked points to produce a degree $15$ curve in $\mathbb{P}^6$. A similar approach yields the unirationality of $\mathcal{H}_{12,8}$, already proven in \cite{SchreyerTanturriMatrix}, and is outlined at the end of Section \ref{unirat137}.

The reversibility of the above constructions corresponds to open conditions on suitable moduli spaces or Hilbert schemes. To show that the so-constructed families of covers of $\mathbb{P}^1$ are dominant on the Hurwitz schemes it is thus sufficient to exhibit single explicit examples of the constructions over a finite field. A computer-aided verification with the computer algebra software \texttt{Macaulay2} \cite{m2} is implemented in the package \cite{package}, whose documentation illustrates the basic commands needed to check the truthfulness of our claims. A ready-to-read compiled execution of our code is also provided.

A priori, it might be possible to mimic these ideas for other pairs $(g,d)$ for which no unirationality result is currently known. However, a case-by-case analysis suggests that, in order to apply the liaison techniques as above, one needs to construct particular curves, which are at the same time far from being general and not easy to realize.

The paper is structured as follows. In Section \ref{preliminaries} we introduce some notation and some background on Brill--Noether Theory and liaison. In Sections \ref{unirat108} and \ref{unirat137} we prove the unirationality of $\mathcal{H}_{10,8}$ and $\mathcal{H}_{13,7}$ respectively.

\begin{ack}
The authors are grateful to Frank-Olaf Schreyer for suggesting the problem and for many useful conversations. The second author also thanks Daniele Agostini for interesting discussions, and the Mathematics and Informatics department of the Universit\"at des Saarlandes for hospitality. The first author is grateful to the DAAD for providing the financial support of her studies. The second author is supported by the Labex CEMPI (ANR-11-LABX-0007-01).
\end{ack}

\section{Preliminaries}
\label{preliminaries}
In this section, we introduce some notation and some background facts which will be needed later on.

\subsection{Brill--Noether Theory}
We recall a few facts from Brill--Noether theory, for which we refer to \cite{ACGH}. Throughout this section, $C$ denotes a smooth general curve of genus $ g $, and $ d,r $ are non-negative integers.

A linear series on $C$ of degree $d$ and dimension $r$, usually referred to as a $g^r_d$, is a pair $(L,V)$, where $L \in \Pic^d(C)$ is a line bundle of degree $d$ and $V \subset \Hb^0(C,L)$ is an $(r+1)$-dimensional vector space of sections of $L$. $C$ has a $g^r_d$ if and only if the Brill--Noether number
\[ \rho=\rho(g,r,d)=g-(r+1)(g+r-d) \]
is non-negative. Moreover, in this case, the Brill-Noether scheme
\[ W^r_d(C)=\{ L \in \Pic^d(C) \mid \hsm^0(L) \ge r+1 \} \]
has dimension $\rho$. The universal Brill--Noether scheme is defined as
\[
\mathcal{W}^r_{g,d} = \{ (C,L) \mid C \in \mathcal{M}_{g}, L \in W^r_d(C) \}.
\]
There is a natural dominant morphism $\alpha:\mathcal{H}_{g,d} \rightarrow \mathcal{W}^1_{g,d}$, which is a $\pgl(2)$-bundle over a dense open subset of $\mathcal{W}^1_{g,d}$; thus, the unirationality of $\mathcal{H}_{g,d}$ is equivalent to the unirationality of $\mathcal{W}^1_{g,d}$, and both $\mathcal{H}_{g,d}$ and $\mathcal{W}^1_{g,d}$ are irreducible.
\subsection{Liaison}
We recall some basic facts on liaison theory.
\begin{definition}
Let $ C $ and $ C^{\prime} $ be two curves in a projective variety $X$ with no embedded and no common components, contained in $r-1$ mutually independent hypersurfaces $Y_i \subset X$ meeting transversally. Let $Y$ be the complete intersection curve $\cap Y_i$. $ C $ and $ C^{\prime} $ are said to be geometrically linked via $Y$ if $ C\cup C^{\prime}=Y$ scheme-theoretically.
\end{definition}
If we assume that the curves are locally complete intersections and that they meet only in ordinary double points, then $\omega_Y|_C =\omega_C(C \cap C')$ and the arithmetic genera of the curves are related by
\begin{equation}
\label{differenceGenera}
2(p_a(C)-p_a(C'))=\deg (\omega_C) - \deg(\omega_{C'})=\omega_X(Y_{1}+\cdots+Y_{r-1}).(C-C').
\end{equation}
The relation above and the obvious relation $\deg C + \deg C' = \deg Y$ can be used to deduce the genus and degree of $C'$ from the genus and degree of $C$.

\medskip
Let $X=\mathbb{P}^1 \times \mathbb{P}^2$ and $C$ be a curve of genus $p_a(C)$ and bidegree $(d_1, d_2)$. With the above hypotheses, let $Y_1, Y_2$ be two hypersurfaces of bidegree $(a_1,b_1)$ and $(a_2,b_2)$, then the genus and the bidegree of $C'$ are
\begin{equation}
\label{liaisonP1P2}
	\begin{array}{rcl}
    	(d_1',d_2') &= &(b_1b_2-d_1,a_1b_2+a_2b_1-d_2),\\
    	p_a(C^{\prime}) &= &p_a(C) - \frac12 \left((a_1+a_2-2)(d_1-d_1')+(b_1+b_2-3)(d_2-d_2')\right).
	\end{array}
\end{equation}

For curves embedded in a projective space $\mathbb{P}^r$, the invariants $p_a(C^{\prime}),d'$ of the curve $C'$ can be computed via
\begin{equation}
\label{liaisonPr}
    \begin{array}{rcl}
    d'&=&\prod d_i - d,\\
    p_a(C^{\prime})&=&p_a(C)-\frac{1}{2}\left(\sum d_i-(r+1)\right)(d-d'),
	\end{array}
\end{equation}
where the $d_i$'s are the degrees of the $r-1$ hypersurfaces $Y_i$ cutting out $Y$.
\section{Unirationality of $\cH_{10,8}$}
\label{unirat108}
In this section we prove the unirationality of $\cH_{10,8}$. To simplify the notation, $\mathbb{P}^1 \times \mathbb{P}^2$ will be denoted by $\PP$. 
\subsection{The double liaison construction}
\label{doubleLiaison}
Let $ (C,L) $ be a general element of $ \cW^{1}_{10,8} $. As $ \rho(10,8,2)< 0 \leq \rho(10,8,1) $, $\hsm^0(L)=2$ and by Riemann--Roch $ \vert K-L\vert $ is a $ 2 $-dimensional linear series of degree $ 10 $. For a general $ 6 $-gonal pencil $ \vert D_1\vert $ of divisors on $ C $, let 

\[ \phi: C\xrightarrow{\vert D_1\vert \times \vert K-L\vert}\PP \]
be the associated map. We assume $ \phi $ is an embedding, and in fact this is the case if the plane model of $ C $ inside $ \mathbb{P}^{2} $ has only ordinary double points and no other singularities, and the points in the preimage of each node under $|K-L|$ are not identified under the map to $\mathbb{P}^{1} $. This way we can identify $ C $ with its image under $ \phi $, a curve of bidegree $(6,10)$ in $\PP$.

Moreover, assume $C$ satisfies the maximal rank condition in bidegrees $(a,3)$ for all $ a\geq 1 $, that is the maps $ \Hb^{0}(\mathcal{O}_{\mathbb{P}}(a,3))\longrightarrow  \Hb^{0}(\mathcal{O}_{C}(a,3))$ are of maximal rank.  Let
$ a_{3} $ be the minimum degree such that $C$ lies on a hypersurface of bidegree $ (a_3,3)$.
Then by Riemann--Roch the maximal rank condition gives $ a_3=3 $ and $C $ is expected to be contained in only one hypersurface of bidegree $(3,3)$. Let $ Y $ be a complete intersection curve containing $ C $ defined by two forms of bidegrees $(3,3)$ and $(4,3)$, and let $ C' $ be the curve linked to $ C $ via $ Y $. By \eqref{liaisonP1P2}, $ C'$ is expected to be a curve of genus $ 4 $ and bidegree $ d^{\prime}=(3,11) $.

Thinking of $ C'$ as a family of three points in $ \mathbb{P}^{2} $ parametrized by the projective line $ \mathbb{P}^{1} $, we expect a finite number $l'$ of distinguished fibers where the three points are collinear. In fact, this is the case when the six planar points of $ C $ lie on a (possibly reducible) conic. We claim that $l'=5$.

To compute $l'$, we need to understand the geometry of $C'$. Let $ D^{\prime}_{2} $ be the divisor of degree $11$ such that the projection of $ C^{\prime} $ to $ \mathbb{P}^{2} $ is defined by a linear subspace of $\Hb^{0}(\mathcal{O}(D^{\prime}_{2}) )$, and let $ \vert D^{\prime}_1\vert $ be the $ 3 $-gonal pencil of divisors defining the map $ C^{\prime}\longrightarrow \mathbb{P}^{1} $. Since $\deg(K_{C'}-D_2')<0$, by Riemann--Roch we have $\hsm^0(\mathcal{O}(D_2'))=11+1-4=8$. We consider the map induced by the complete linear system
\begin{equation}
\label{definepsi2}
\psi_{2}:C^{\prime}\xrightarrow{|D_2'|} \mathbb{P}^{7};
\end{equation}
as shown in \cite{Sh86}, the $ 3 $-dimensional rational normal scroll $ S $ of degree $ 5 $ swept out by $ \vert D^{\prime}_1\vert $ contains the image of $ \psi_{2} $. Hence, the image of the map
\[ \psi:C^{\prime}\longrightarrow \mathbb{P}^{1}\times \mathbb{P}^{7} \]
is contained in the graph of the natural projection map from $ S$ to $ \mathbb{P}^{1} $, that is 
\[ \psi(C^{\prime})\subseteq \mathbb{P}^{1}\times S=\bigcup_{D_\lambda \in \vert D^{\prime}_1\vert}([\lambda]\times \bar{D}_\lambda ), \]
where $ \bar{D}_\lambda $ is the linear span of $ \psi_{2}( D_\lambda) $ in $ \mathbb{P}^{7} $.

As $ \psi(C^{\prime}) $ is a family of three points in $ \mathbb{P}^{7} $ parametrized by $ \mathbb{P}^{1} $, $ C^{\prime}\subset \mathbb{P}^{1}\times \mathbb{P}^{2} $ is obtained by projection of $ \psi(C^{\prime}) $ from a linear subspace $\mathbb{P}^{1}\times V \subset \mathbb{P}^{1}\times \mathbb{P}^{7}$ of codimension $ 3 $. Fix a $\lambda \in \mathbb{P}^1$; by Riemann--Roch, $\dim |D_2'-D_1'|=5$, hence $\psi_2(D_\lambda)$ spans a 2-dimensional projective space inside $\mathbb{P}^7$. It is clear that the three points corresponding to $\lambda$ are distinct and collinear if and only if $V \cap \bar{D_\lambda}$ is a point, and the three points coincide if and only if $V \cap \bar{D_\lambda}$ is a projective line. The latter case does not occur in general, as the plane model of $C'$ has only double points. The former case occurs in $l'=\deg S=5$ points if $S$ and $V$ intersects transversally, an open condition which holds in general.

Now, suppose that for all $ b\geq 1$ the maps 
\[ \Hb^{0}(\mathcal{O}_{\mathbb{P}}(b,2))\longrightarrow  \Hb^{0}(\mathcal{O}_{C^{\prime}}(b,2)) \]
are of maximal rank, and set 
\[ b_{2}:=\min\lbrace b\ : \ \hsm^{0}(\mathcal{I}_{C'}(b,2))\neq 0\rbrace.\]
Under the maximal rank assumption, $ b_{2}=5 $ and $ \hsm^{0}(\mathcal{I}_{C'}(5,2))=2 $. Let $ Y'$ be a complete intersection of two hypersurfaces of bidegree $ (5,2)$ containing $C'$, and let $ C^{\prime\prime} $ be the curves linked to $ C'$ via $ Y'$.

Interpreting again $C'$ and $C''$ as families of points parametrized by $\mathbb{P}^1$, we observe that a general fiber of $ C'' $ consists of a single point. In the $5$ distinguished fibers of $C'$, the two conics of the complete intersection $ Y'$ turn out to be reducible with the line spanned by the three points of $ C' $ as common factor. Thus, the curve $ C'' $ is the union of a rational curve $ R $ of bidegree $ (1,4) $ and $5 $ lines.

\subsection{A unirational parametrization}
The double liaison construction described above can be reversed and implemented in a computer algebra system. We note that all the assumptions on $C$ and $C'$ we made correspond to open conditions in suitable moduli spaces or Hilbert schemes, so that it is sufficient to check them on a single example. We can work on a finite field, as explained in Remark \ref{finiteField} here below.

\begin{remark}
	\label{finiteField}
In this paper we will often need to exhibit an explicit example satisfying some open conditions. A priori we could perform our computations directly on $\mathbb{Q}$, but this can increase dramatically the required time of execution. Instead, we can view our choice of the initial parameters in a finite field $F_p$ as the reduction modulo $p$ of some choices of parameters in $\mathbb{Z}$. Then, the so-obtained example $E_p$ can be seen as the reduction modulo $p$ of a family of examples defined over a neighborhood $\Spec \mathbb{Z}[\frac{1}{b}]$ of $(p) \in \Spec \mathbb{Z}$ for a suitable $b \in \mathbb{Z}$ with $ p \nmid  b$. If our example $E_p$ satisfies some open conditions, then by semicontinuity the generic fiber $E$ satisfies the same open conditions, and so does the general element of the family over $\mathbb{Q}$ or $\mathbb{C}$.
\end{remark}

Our construction depends on a suitable number of free parameters corresponding to the choices we made. Picking $ 5 $ lines in $ \mathbb{P}^{1}\times \mathbb{P}^{2} $ requires $ 5\cdot 3=15 $ parameters. Choosing $ 2 $ forms of bidegree $ (2,1) $ to define the rational curve $R$ corresponds to the choice of $\dim \Gr(2,9) = 14$ parameters. By Riemann--Roch we expect $  \hsm^{0}(\mathcal{I}_{C''}(5,2))=7$, so we need $ \dim \Gr(2,7)=10 $ parameters to define the complete intersection $ Y'$. Similarly, as $ \hsm^{0}(\mathcal{I}_{C'}(3,3))=1 $ and $ \hsm^{0}(\mathcal{I}_{C'}(4,3))=8 $, we require $ \dim \Gr(1,8)-2=5 $ further parameters for the complete intersection $ Y $. This amounts to $ 15+14+10+5=44 $ parameters in total.

\begin{theorem}
	\label{H108unirational}
The Hurwitz space $\cH_{10,8}$ is unirational.
\begin{proof}
Let $\mathbb{A}^{44} $ be our parameter space. With the code provided by the function \texttt{verifyAssertionsOfThePaper(1)} in \cite{package}, and following the construction of Section \ref{doubleLiaison} backwards we are able to produce an example of a curve $C \subset \PP$ and to check that all the assumptions we made are satisfied, that is:
\begin{itemize}
\item for a general choice of a curve $C''$, a union of a rational curve of bidegree $(1,4)$ and $5$ lines, and for a general choice of two hypersurfaces of bidegree $(5,2)$ containing $C''$, the residual curve $C'$ is a smooth curve of genus $4$ and bidegree $(3,11)$ which intersects $C''$ only in ordinary double points;
\item $C'$ satisfies the maximal rank condition in bidegrees $(b,2)$ for all $b \geq 1$ and its planar model has only ordinary double points as singularities;
\item for a general choice of two hypersurfaces of bidegree $(3,3), (4,3)$ containing $C'$, the residual curve $C$ is a smooth curve of genus $10$ and bidegree $(6,10)$ that intersects $C'$ only in ordinary double points;
\item $C$ satisfies the maximal rank condition in bidegrees $(a,3)$ for all $a \geq 1$ and its planar model is non-degenerate.
\end{itemize}
This means that our construction produces a rational family of elements in $\mathcal{W}^2_{10,10}$, the Serre dual space to $\mathcal{W}^1_{10,8}$. As all the above conditions are open, and  $\mathcal{W}^1_{10,8}$ is irreducible, this family is dominant, which proves the unirationality of both $\mathcal{W}^1_{10,8}$ and $\mathcal{H}_{10,8}$. \qedhere

\end{proof}
\end{theorem}

\section{Unirationality of $\cH_{13,7}$}
\label{unirat137}
In this section we will prove the unirationality of the Hurwitz space $\cH_{13,7}$. As a preliminary result of independent interest, let us prove the following

\begin{theorem}
	\label{M10n}
	The moduli space $\mathcal{M}_{10,n}$ of curves of genus $10$ with $n$ marked points is unirational for $1 \leq n \leq 5$. 
	\begin{proof}
	This result is achieved by linkage on $\mathbb{P}:=\mathbb{P}^1 \times \mathbb{P}^2$. We start with a reducible curve $C$ of arithmetic genus $-3$, union of 3 general lines and the graph of a rational plane curve of degree $4$. On the one hand, the space of such curves is clearly unirational; on the other hand, in general $C$ will be contained in at least two independent hypersurfaces of bidegree $(4,2)$. The linkage with respect to 2 general such hypersurfaces produces a curve $C'$ of expected bidegree $(3,9)$ and genus $4$, which will be in general contained in exactly $7$ independent hypersurfaces of bidegree $(3,3)$.
	
	For the choice of 5 general points $\{P_1,\dotsc,P_5\}$ in $\mathbb{P}$, let $I_P$ be their ideal. In general, the space of bihomogeneous polynomials $(3,3)$ contained in $I_{C'} \cap I_P$ will be generated by two independent polynomials $f_1,f_2$, defining two hypersurfaces $X_1,X_2$. The complete intersection of these hypersurfaces link $C'$ to a curve $C''$ passing through each $P_i$; $C''$ turns out to be a curve of genus $10$ and bidegree $(6,9)$. The projection onto $\mathbb{P}^1$ yields an element of $\cH_{10,6}$.
	
	In \cite{GeissUnirationality} Gei\ss\ proved that this construction yields a rational dominant family in $\cH_{10,6}$. Moreover, the Brill--Noether number $\rho(6, 1, 10)= 10 - (1+1)(10 - 6+1)=0$ is non-negative, which implies that this rational family dominates $\mathcal{M}_{10}$ as well. Therefore, as the choice of $\{P_1,\dotsc,P_5\}$ is unirational, we get a rational dominant family of curves of genus $10$ together with (up to) five marked points.
	\end{proof}
\end{theorem}

\begin{theorem}
	\label{H137unirational}
The Hurwitz space $\cH_{13,7}$ is unirational.
\begin{proof}
	Let $(D,L) \in \cW^1_{13,7}$ be a general element. By Riemann--Roch, $\omega_{D}\otimes L^{-1}$ is a general $g^6_{13,17}$ and therefore the linear system $|K_D-L|$ embeds $D$ in $\mathbb{P}^6$ as a curve of genus $13$ and degree $17$. Conversely, if $D$ is a general curve of genus $13$ and degree $17$ in $\mathbb{P}^6$, by Riemann--Roch the line bundle $\omega_{D}\otimes\mathcal{O}_{D}(-1)$ is a general $g^1_7$. Hence, in order to prove the unirationality of $\cH_{13,7}$, it will be sufficient to exhibit a rational family of projective curves of genus $13$ and degree $17$ in $\mathbb{P}^6$ which dominates $\cW^6_{13,17}$.
	
	Let $C$ be a general curve of genus $13$ and degree $17$ in $\mathbb{P}^6$. Since $\mathcal{O}_C(2)$ is non-special, $C$ is contained in at least ${{6+2}\choose{2}}-(17\cdot 2 +1 - 13)=6$ independent quadric hypersurfaces. Consider five general such hypersurfaces $X_i$ and suppose that the residual curve $C'$ is smooth and that $C$ and $C'$ intersect transversally; these are open conditions on the choice of $(C,\mathcal{O}_C(1)) \in \cW^6_{13,17}$. By \eqref{liaisonPr}, $C'$ has genus $g'=10$ and degree $d'=15$.
	By Riemann--Roch, the Serre residual divisor $\omega_{C'}\otimes\mathcal{O}_{C'}(-1)$ has degree 3 and one-dimensional space of global sections, hence it corresponds to the class of three points on $C'$. Conversely, by Geometric Riemann--Roch three general points on $C'$ form a divisor $P$ with $\hsm^0(P)=1$ such that $|K_C-P|$ embeds $C'$ in $\mathbb{P}^6$ as a curve of degree 15. Hence, the unirationality of $\mathcal{W}^6_{10,15}$ can be deduced from the unirationality of $\mathcal{M}_{10,3}$, proved in Theorem \ref{M10n} above.
	
	By means of the implemented code \texttt{verifyAssertionsOfThePaper(2)} in \cite{package}, we can show with an explicit example that
	\begin{itemize}
	\item for a general curve $C'$ of genus 10 and degree 15 in $\mathbb{P}^6$ and for a general choice of five quadric hypersurfaces containing it, the residual curve $C$ is smooth and intersects $C'$ only in ordinary double points;
	\item $C$ is not contained in any hyperplane.
	\end{itemize}
	This way we get a rational family of curves $C$ of genus 13 and degree 17 in $\mathbb{P}^6$. Since all the assumptions we made correspond to open conditions on $\mathcal{W}^6_{13,17}$ and are satisfied by our explicit examples, such family dominates $\mathcal{W}^6_{13,17}$.
\end{proof}
\end{theorem}

\begin{remark}
The same argument of Theorem \ref{H137unirational} holds for a general element in $\cH_{12,8}$, so that the above proof yields an alternative proof of the unirationality of $\cH_{12,8}$ proved in \cite{SchreyerTanturriMatrix}. In this case, the Serre dual model is a curve of genus $12$ and degree $14$ in $\mathbb{P}^4$. The liaison is taken with respect to 3 general cubic hypersurfaces and yields a curve of genus 10 and degree 13, which can be constructed from a curve of genus 10 and 5 marked points with the same strategy as above. An implementation of this unirational parametrization of $\cH_{12,8}$ via linkage can be found in the package \cite{package}.
\end{remark}

The package \cite{package} including the implementation of the unirational pa\-ra\-me\-tri\-za\-tions exhibited in the paper, together with all the necessary and supporting documentation, is available online.

%\bibliographystyle{alpha}
%\bibliography{database}

%\makeatletter
%\providecommand\@dotsep{5}
%\makeatother
%\listoftodos\relax

\end{document}